\documentclass{ifacconf}
\usepackage{graphicx}      
\usepackage{natbib}        
\usepackage{amsfonts}
\usepackage{algorithmic}
\usepackage{amsmath,amssymb} 
\usepackage{url,cite}
\usepackage{graphicx}
\usepackage{tikz,pgfplots}
\usepackage{mathrsfs}

\allowbreak
\allowdisplaybreaks

\usepackage{listings}

\lstdefinelanguage{julia}{
    alsoletter={0123456789!.:+-^?*~><='/},
    keywords={},
    keywordsprefix=:,
    keywords = [2]{function, end, begin, if, elseif, else, using, for},
    keywords = [3]{InfiniteModel, MvNormal, optimize!, objective_value, Infinite, zeros,
                   DomainRestrictions, set_value, enumerate, value, value., length, supports,
                   OrthogonalCollocation, sum, add_supports, Uniform},
    keywords = [4]{0, 1, 0.303, 0.727, 0.3, 1e-5, 200, 101, 0.02, 0.8},
    keywords = [5]{@variable, @constraint, @objective, @parameter_function, 
                   @infinite_parameter, @finite_parameter, @variables, @constraints},
    keywords = [6]{+, -, ^, ?, *, >, <, =, <=, >=, ==, .<=, .>=, .==, .+, .-, =>, in, /},
    keywordsprefix=:,
    keywordstyle=\color{red},
    keywordstyle=[2]\color{violet},
    keywordstyle=[3]\color{blue},
    keywordstyle=[4]\color{orange},
    keywordstyle=[5]\color{purple},
    keywordstyle=[6]\color{teal},
    sensitive=true,
    morecomment=[l]{\#},
    morestring=[b]"
}

\begin{document}
\begin{frontmatter}

\title{New Measures for Shaping Trajectories \\in Dynamic Optimization}  

\thanks[footnoteinfo]{This work was supported by the U.S. Department of Energy under grant DE-SC0014114.}

\author{Joshua L. Pulsipher, Benjamin R. Davidson, and Victor M. Zavala}

\address{Department of Chemical and Biological Engineering,\\
  University of Wisconsin-Madison, Madison, WI 53706 USA \\
  (e-mail: \{pulsipher,victor.zavala\}@wisc.edu).}

\begin{abstract} 
We propose a new class of measures for shaping time-dependent trajectories in dynamic optimization (DO). The proposed measures are analogous to risk measures used in stochastic optimization (SO) and are inspired by a recently-proposed unifying abstraction for infinite-dimensional optimization. Risk measures are summarizing statistics (e.g., average, variance, quantiles, worst-case values) that are used to shape the probability density of random objectives and constraints. We show that this extensive collection of measures can be applied in DO for computing and manipulating  interesting features of time-dependent trajectories (e.g., excursion costs and quantiles). We also discuss how to implement these measures in the Julia modeling package \texttt{InfiniteOpt.jl}. 
\end{abstract}


\end{frontmatter}

\section{Introduction} \label{sec:intro}
We consider the continuous-time dynamic optimization (DO) problem:
\begin{equation}
    \begin{aligned}
    && \min_{y(\cdot) \in \mathcal{Y}} &&& M_t f(\dot{y}(t), y(t), t) \\
    && \text{s.t.} &&& g(\dot{y}(t), y(t), t) \leq 0, && t \in \mathcal{D}_t. \\
    \end{aligned}
    \label{eq:gen_form}
\end{equation}
This representation follows from a recently-proposed abstraction that unifies infinite-dimensional optimization (InfiniteOpt) problems (e.g., dynamic/stochastic/PDE optimization) \citep{pulsipher2021unifying}. In the context of DO, the infinite domain $\mathcal{D}_t$ is the time domain $[t_0, t_f]$, $y(t) \in \mathcal{Y} \subseteq \mathbb{R}^{n_y}$ are time-valued decision functions (e.g., state/control variables), $\dot{y}(t) \in \mathbb{R}^{n_y}$ are derivative variables $dy(t)/dt$, $f(\cdot)$ is an infinite-dimensional cost function, and $g(\cdot)$ is a vector-valued infinite-dimensional constraint functions $g_j(\cdot), \; j \in \mathcal{J} \subseteq \mathbb{R}^{n_g}$. 

The measure operator $M_t : \mathcal{D}_t \mapsto \mathbb{R}$ is the focus of this work; this operator seeks to scalarize infinite-dimensional cost functions (summarizing them over the time domain $\mathcal{D}_t$) to form a well-posed objective function. Such measures can also be used to handle constraints, but here we focus on objectives to simplify the presentation. Problem \eqref{eq:gen_form} is general and captures a wide range of DO problems (e.g., model predictive control and state/parameter estimation) where $y(\cdot)$ is comprised of state/control/parameter variables and $g(\cdot) \leq 0$ can include DAE/path/point constraints. 

Classical DO formulations minimize the integral of the cost trajectory (known as the Bolza objective). In other words, these formulations use the measure operator:
\begin{equation}
    M_t f(t)=\int_{t \in \mathcal{D}_t} f(t) dt
    \label{eq:trad_obj}
\end{equation}
where we write $f(t) := f(\dot{y}(t), y(t), t)$ for convenience. Minimizing \eqref{eq:trad_obj} has the effect of {\em uniformly} shaping the cost trajectory $f(t)$ over the domain $\mathcal{D}_t$. This also amounts to minimizing the total cost (which is equivalent to the average cost with a normalization factor of $(t_f-t_0)^{-1}$). However, we can envision more advanced DO formulations that aim to shape the cost trajectory $f(\cdot)$ with other measures (e.g., peak or excursion costs). For example, Risbeck and Rawlings recently proposed an MPC objective that penalizes the total and peak of the cost \citep{risbeck2019economic}. Minimizing the peak of a cost-trajectory is often a desirable feature in systems, as these can be associated with unsafe behavior and economic penalties. 

Under the InfiniteOpt representation, Problem \eqref{eq:gen_form} is analogous to the SO problem:
\begin{equation}
    \begin{aligned}
    && \min_{z \in \mathcal{Z}, y(\cdot) \in \mathcal{Y}} &&& M_\xi f(z, y(\xi), \xi) \\
    && \text{s.t.} &&& g(z, y(\xi), \xi) \leq 0, && \xi \in \mathcal{D}_\xi. \\
    \end{aligned}
    \label{eq:2stage}
\end{equation}
where $z \in \mathcal{Z} \subseteq \mathbb{R}^{n_z}$ are here-and-now decision variables, $\xi$ is a random parameter, $\mathcal{D}_\xi \subseteq \mathbb{R}^{n_\xi}$ is the co-domain of its distribution, and $y(\xi) \in \mathcal{Y} \in \mathbb{R}^{n_y}$ are recourse decision variables. We observe that \eqref{eq:gen_form} is a special case of \eqref{eq:2stage} if we define $n_z = 0$, $n_\xi = 1$, $\mathcal{D}_\xi = \mathcal{D}_t$, $M_\xi = M_t$, and set $y(\xi)$ to contain both $\dot{y}(t)$ and $y(t)$. 

In SO, the measure $M_\xi$ can draw from a wide collection of risk measures (summarizing statistics) that aim to shape the probability density function of $f(\xi)$; some examples include the expected value, mean-variance, absolute deviation, and conditional-value-at-risk (CVaR) \citep{ruszczynski2006optimization}. Risk measures enable greater flexibility in shaping the cost function $f(\xi)$ and are commonly used to penalize extreme events (i.e., high costs). By leveraging the connection between \eqref{eq:gen_form} and \eqref{eq:2stage}, we have recently shown that the expected value and CVaR measures can be used as measures in DO to shape the time-dependent cost trajectory  \citep{pulsipher2021unifying}. 

In this work, we formalize the use of risk measures to DO problems that follow \eqref{eq:gen_form}. This new class of dynamic measures $M_t$ enable us to shape trajectories and policies in new ways. Moreover, we highlight how the Julia-based modeling package \texttt{InfiniteOpt.jl} provides a useful interface to compactly express such measures. 

\vspace{-0.05in}

\section{Main Results} \label{sec:main}

\subsection{Time-Valued Density Functions}

We beging by establishing time analogues of probability density functions (pdf) and of cumulative density functions (cdf) typically used in SO (i.e., pdfs and cdfs defined over the time domain $\mathcal{D}_t$). These provide key constructs to facilitate the use of risk measures in DO.  

In SO, $\xi$ is described by the pdf $p_\xi : \mathcal{D}_\xi \mapsto \mathbb{R}_{\geq0}$, which satisfies $\int_{\xi \in \mathcal{D}_\xi} p_\xi(\xi) d\xi = 1$. The pdf is used to compute statistics such as the expectation:
\begin{equation}
    \mathbb{E}_\xi[f(\xi)] = \int_{\xi \in \mathcal{D}_\xi} f(\xi) p_\xi(\xi) d\xi.
    \label{eq:random_expect}
\end{equation}
This measure $M_\xi = \mathbb{E}_\xi$ is commonly used in SO; here, we can see that $p_\xi(\cdot)$ acts as a weighting function that places varied emphasis over the domain $\mathcal{D}_\xi$. The expectation summarizes the cost function in a single scalar value. 

Following the analogy of SO and DO, we define a weighting function $p_t : \mathcal{D}_t \mapsto \mathbb{R}_{\geq0}$ (i.e., a time-valued pdf) to yield the time average (expectation):
\begin{equation}
    \mathbb{E}_t[f(t)] := \int_{t \in \mathcal{D}_t} f(t) p_t(t) dt.
    \label{eq:time_expect}
\end{equation}
This provides us flexibility in prioritizing different regimes in the domain $\mathcal{D}_t$. We obtain  \eqref{eq:trad_obj} as a special case by setting $p(t) = 1$; however, we can envision choosing from a wide range of candidate pdfs, such as:
\begin{equation}
    p_t(t) = \frac{1}{t_f - t_0}
    \label{eq:uniform_pdf}
\end{equation}
Note that this pdf of the uniform random parameter and places equal emphasis on different parts of the time domain. This inspires the consideration of other weighting functions, such as the exponential pdf:
\begin{equation}
    p_t(t) = \gamma \e^{-\gamma t}
    \label{eq:exp_pdf}
\end{equation}
where $\gamma \in \mathbb{R}_{>0}$ is the decay rate parameter.  Using this weighting function to compute the time average places emphasis on the beginning of the domain and quickly decays over time. In other words, this can be seen as a discount factor \citep{shin2020diffusing}. 

The cdf $P : \mathcal{D}_\xi \mapsto [0, 1]$  of a random variable $\xi$ is:
\begin{equation}
    P(\xi; \hat{\xi}) := \mathbb{P}_\xi(\xi \leq \hat{\xi}) = \int_{\xi \in \{\xi \in \mathcal{D}_\xi : \xi \leq \hat{\xi} \}} p_\xi(\xi) d\xi
    \label{eq:random_cdf}
\end{equation}
This denotes the cumulative probability of finding $\xi$ below the threshold $\hat{\xi} \in \mathcal{D}_\xi$. In the context of the DO, we can write the cdf of the time trajectory $f(t)$ by using the {\em excursion sets}: 
\begin{equation}
    \begin{gathered}
        \mathcal{D}_t^+(f(t); \hat{f}) := \{t \in \mathcal{D}_t : f(t) \geq \hat{f}\} \\
        \mathcal{D}_t^-(f(t); \hat{f}) := \{t \in \mathcal{D}_t : f(t) \leq \hat{f}\}
    \end{gathered}
\end{equation}
where $\mathcal{D}_t^+(f(t); \hat{f}) \subseteq \mathcal{D}_t$ and $\mathcal{D}_t^-(f(t); \hat{f}) \subseteq \mathcal{D}_t$ are the positive and negative function excursion sets, respectively. We can use the negative excursion set to establish:
\begin{equation}
    P(f(t); \hat{f}) = \int_{t \in \mathcal{D}_t^-(f(t); \hat{f})} p_t(t) dt.
    \label{eq:time_cdf}
\end{equation}
In a DO context, the cdf measures the fraction of time that the trajectory $f(t)$ is below the threshold $\hat{f}$. We will  see that expressing the cdf in terms of excursion set allows us to interpret measures as mechanisms to bound  trajectories.  

\subsection{Dynamic Measures} \label{sec:measures}
Here we illustrate the interpretation of risk measures in a DO context. A large number of risk measures have been proposed in the SO community and analyzing them is beyond the scope of this work (we refer the reader to \citep{krokhmal2013modeling} for a review). This section establishes key constructs and exemplifies the steps needed to interpret risk measures in a time setting. We discuss properties of these proposed measures in Section \ref{sec:properties}. 

\subsubsection{Expectation}

The time expectation measure $\mathbb{E}_t$ shown in Equation \eqref{eq:time_expect} allows us to assess a weighted average of our cost trajectory $f(t)$ in accordance with the weighting function $p_t(t)$. In the context of DO, minimizing the expected cost provides an easily interpretable objective and the key modeling choice lies in the selection of $p_t(t)$.  Reasonable candidates for many DO applications are \eqref{eq:uniform_pdf} and \eqref{eq:exp_pdf}, as demonstrated in Section \ref{sec:cases}; however, other weighting functions are possible (e.g., Gaussian and Gamma). The selection of the weighting function dictates how much emphasis is placed on different parts of the time domain. The expectation serves as a core construct in defining the more sophisticated measures. 

\begin{prop} \label{prop:expect_equiv}
The time expectation $\mathbb{E}_t$ in \eqref{eq:time_expect} is a special case of $\mathbb{E}_\xi$ in  \eqref{eq:random_expect} if $n_z = 0$, $n_\xi = 1$, $\mathcal{D}_\xi = \mathcal{D}_t$. 
\end{prop}

\subsubsection{Mean-Variance} 

The mean-variance $\mathbb{E}$-$\mathbb{V}_\xi$ is a classical measure used in portfolio optimization \citep{leland1999beyond}. For SO problems  \eqref{eq:2stage}, this minimizes the variance (i.e., spread) of the cost outcomes in an attempt to mitigate high-cost events: 
\begin{equation}
    \mathbb{E}\text{-}\mathbb{V}_\xi[f(\xi)] := \mathbb{E}_\xi[f(\xi)] + \lambda \mathbb{V}_\xi[f(\xi)]
    \label{eq:random_mean_var}
\end{equation}
where $\mathbb{V}_\xi = \mathbb{E}_\xi[(f(\xi) - \mathbb{E}_\xi[f(\xi)])^2]$ is the variance and $\lambda \in \mathbb{R}_{\geq0}$ is a tradeoff parameter. Transferring this to a DO setting we obtain the time-valued measure:
\begin{equation}
    \mathbb{E}\text{-}\mathbb{V}_t[f(t)] := \mathbb{E}_t[f(t)] + \lambda \mathbb{V}_t[f(t)].
    \label{eq:time_mean_var}
\end{equation}
Minimizing $\mathbb{E}\text{-}\mathbb{V}_t[f(t)]$ in a DO problem provides a tradeoff problem that seeks to minimize the magnitude of the cost trajectory (the expectation) and the variability/fluctuations (variance) of the cost trajectory. A disadvantage of this measure is that it penalizes cost variability equally for low and high costs. However, this property can be advantageous for cost functions that seek to enforce smooth control trajectories.

\begin{prop} \label{prop:mean_var_equiv}
The measure operator $\mathbb{E}$-$\mathbb{V}_t$ from \eqref{eq:time_mean_var} is a special case of $\mathbb{E}$-$\mathbb{V}_\xi$ from \eqref{eq:random_mean_var} under the same conditions of Proposition \ref{prop:expect_equiv}.
\end{prop}

\subsubsection{Quantile}

The quantile $Q_\xi(f(\xi); \alpha)$ (also referred to as the value-at-risk) denotes the threshold value $\hat{f}$ for $f(\xi)$ such that the cumulative probability of incurring costs below the threshold is at least $\alpha \in [0, 1]$:
\begin{equation}
    Q_\xi(f(\xi); \alpha) := \inf_{\hat{f} \in \mathbb{R}} \left\{P(f(\xi); \hat{f}) \geq \alpha \right\}.
    \label{eq:random_quantile}
\end{equation}
Constraining \eqref{eq:random_quantile} is equivalent to enforcing a probabilistic constraint  \citep{sarykalin2008value}:
\begin{equation}
    Q_\xi(f(\xi); \alpha) \leq 0 \iff \mathbb{P}_\xi(f(\xi) \leq 0) \geq \alpha.
\end{equation}
We can use the cdf \eqref{eq:time_cdf} in combination with \eqref{eq:random_quantile} to define the time-valued quantile:
\begin{equation}
    Q_t(f(t); \alpha) := \inf_{\hat{f} \in \mathbb{R}} \left\{\int_{t \in \mathcal{D}_t^-(f(t); \hat{f})} p_t(t) dt \geq \alpha \right\}.
    \label{eq:time_quantile}
\end{equation}
Using this measure in \eqref{eq:gen_form} minimizes the excursion threshold of the cost function trajectory such that the fraction of time that exceed it is no more than $1 - \alpha$. Unlike the $\mathbb{E}$-$\mathbb{V}_t$ measure, the quantile measure only penalizes high cost values, making it an attractive alternative in certain cases. However, the potential disadvantages of this measure are that it does not strongly discourage high cost peaks in the positive function excursion $\mathcal{D}_t^+(f(t); Q_t(f(t); \alpha))$ and it is nonconvex and difficult to compute in general. 

\begin{prop} \label{prop:quantile_equiv}
The quantile measure $Q_t$ \eqref{eq:time_quantile} is a special case of its analogue $Q_\xi$ in  \eqref{eq:random_quantile} under the same conditions of Proposition \ref{prop:expect_equiv}.
\end{prop}

\subsubsection{Conditional-Value-at-Risk}
The conditional-value-at-risk (CVaR) measure seeks to address the limitations of the quantile measure $\mathbb{E}$-$\mathbb{V}_\xi$ by penalizing the expected value of the $1-\alpha$ largest cost values:
\begin{equation}
    \text{CVaR}_\xi(f(\xi); \alpha) := \min_{\hat{f} \in \mathbb{R}} \left\{\hat{f} + \frac{1}{1-\alpha} \mathbb{E}_\xi[f(\xi) - \hat{f}]_+ \right\}
    \label{eq:random_cvar}
\end{equation}
where $\mathbb{E}_\xi[f(\xi) - \hat{f}]_+ := \mathbb{E}_\xi[\max(0, f(\xi) - \hat{f})]$ and $\alpha \in [0, 1)$. CVaR is also known as the superquantile; under mild assumptions, CVaR can be represented as:
\begin{equation}
    \text{CVaR}_\xi(f(\xi); \alpha) = \mathbb{E}_\xi[f(\xi) : f(\xi) \geq Q_t(f(t); \alpha)]
    \label{eq:random_cvar2}
\end{equation}
since the minimizer $\hat{f}^*$ is $Q_t(f(t); \alpha)$ \citep{rockafellar2000optimization}. With this observation, the time-valued CVaR measure can be expressed as:
\begin{equation}
    \text{CVaR}_t(f(t); \alpha) := \min_{\hat{f} \in \mathbb{R}} \left\{\hat{f} + \frac{1}{1-\alpha} \mathbb{E}_t[f(t) - \hat{f}]_+ \right\}.
    \label{eq:time_cvar}
\end{equation}
This provides a convex measure that penalizes the high (peak) costs incurred in the positive function excursion set $\mathcal{D}_t^+(f(t); Q_t(f(t); \alpha))$. Note that this penalizes {\em multiple} peak costs (and not just the peak cost, as done in typical DO formulation). Moreover, one can show that:
\begin{equation}
    \begin{aligned}
        &\lim_{\alpha \rightarrow 0} \text{CVaR}_t(f(t);\alpha) = \mathbb{E}_t[f(t)] \\
        &\lim_{\alpha \rightarrow 1} \text{CVaR}_t(f(t);\alpha) = \max_{t \in \mathcal{D}_t} f(t)
    \end{aligned}
    \label{eq:cvar_special}
\end{equation}
which both follow from Equation \eqref{eq:random_cvar2}. As such, CVaR is highly versatile measure for use in DO that can capture both average and extreme features of a time trajectory. 

\begin{prop} \label{prop:cvar_equiv}
The CVaR measure in Equation \eqref{eq:time_cvar} is a special case of CVaR$_\xi$ in Equation \eqref{eq:random_cvar} under the same conditions of Proposition \ref{prop:expect_equiv}.
\end{prop}

\subsubsection{Disutility}
Disutility risk measures are another prevalent measure class used in SO; this employs an expectation over a disutility function $g : \mathbb{R} \mapsto \mathbb{R}$ (typically a convex increasing function) that penalizes unfavorable values of $f(\xi)$:
\begin{equation}
    D_\xi(f(\xi)) := \mathbb{E}_\xi[g(f(\xi))].
    \label{eq:random_disutility}
\end{equation}
In an effort to make a translation invariant measure, the measure is often expressed as:
\begin{equation}
    \tilde{D}_\xi(f(\xi)) := \inf_{\hat{f} \in \mathbb{R}}\mathbb{E}_\xi[f(\xi) + g(f(\xi) - \hat{f})].
    \label{eq:random_disutility2}
\end{equation}
One can show that $\text{CVaR}_\xi(f(\xi); \alpha)$ is a special case of \eqref{eq:random_disutility2}. This follows by 
letting $g(x; \alpha) = (1 - \alpha)^{-1}\max(0, x) - x$ where $x \in \mathbb{R}$. Then by substituting $g(x; \alpha)$ in  \eqref{eq:random_disutility2} we obtain \eqref{eq:random_cvar}. We transfer \eqref{eq:random_disutility2} to DO by using the time-valued expectation $\mathbb{E}_t$:
\begin{equation}
    \tilde{D}_t(f(t)) := \inf_{\hat{f} \in \mathbb{R}}\mathbb{E}_t[f(t) + g(f(t) - \hat{f})].
    \label{eq:time_disutility}
\end{equation}
This measure class provides great flexibility in shaping dynamic trajectories as there are diverse choices of $g(\cdot)$ (in addition to the flexibility provided via selecting the weighting function $p_t(\cdot)$). A useful survey on the properties of disutility functions in the context of SO is provided in \citep{fulga2016portfolio}.  

\begin{prop} \label{prop:disutility_equiv}
The time-valued disutility measure $\tilde{D}_t(f(t))$ is a special case of $\tilde{D}_\xi(f(\xi))$ under the same conditions of Proposition \ref{prop:expect_equiv}.
\end{prop}

\subsection{Measure Properties} \label{sec:properties}
Here we formalize some key mathematical properties of the dynamic measures presented in Section \ref{sec:measures}. This provides some interesting and useful insights on the behavior that these measures induce. These properties have been studied in the SO community and we will show that this rich theory can be readily applied to DO. 

In the context of SO, four main properties are typically considered for risk measures: convexity, monotonicity, translation invariance, and positive homogeneity. Moreover, a measure operator is said to be coherent if it satisfies all these  properties \citep{artzner1999coherent, ruszczynski2006optimization}.  

Convexity asserts that a measure operator $M_\xi$ satisfy:
\begin{equation}
    M_\xi(\beta f + (1 - \beta) h) \leq \beta M_\xi(f) + (1 - \beta) M_\xi(h)
    \label{eq:basic_convexity}
\end{equation}
for all measurable functions $f(\xi),h(\xi) : \mathcal{D}_\xi \mapsto \mathbb{R}$ in the linear function space $\mathscr{F}$ and all $\beta \in [0, 1]$. This property is key for creating optimization objectives that are well-posed and guarantees that the measure of a convex cost is also convex. 

Under monotonicity, we have that if $f_1(\xi)\succeq f_2(\xi)$ ($f_1(\xi)$ dominates $f_2(\xi)$), then the measure $M_\xi$ satisfies:
\begin{equation}
    M_\xi(f_1(\xi)) \geq M_\xi(f_2(\xi)).
    \label{eq:basic_monotonicity}
\end{equation}
This ensures that, if a cost function dominates another cost function, then the measure former will also be greater. than that of the latter. The concept of dominance (comparing whether a random variable is better than another random variable) is an interesting and important concept that has not been explored in DO. In a DO context, dominance of first-order ($f_1(t)\succeq f_2(t)$) requires that $P(f_1(t)>\hat{f})\geq P(f_2(t)>\hat{f})$ for any threshold value $\hat{f}$. In other words, the fraction of time that the trajectory $f_1(t)$ remains above the threshold is greater or equal than the fraction of time that the trajectory $f_2(t)$ remains above the same threshold. A monotonic measure is such that, if $f_1(t)\succeq f_2(t)$, then $M_t(f_1(t)) \geq M_t(f_2(t))$ holds. Note that dominance holds trivially if $f_1(t)\geq f_2(t)$ for all $t\in \mathcal{D}_t$.  These concepts are important because comparisons (benchmarks) of time trajectories are not as straightforward \citep{renteria2018optimal}, as the trajectories are functions (not scalar values) and thus a trajectory might be better in some parts of the time domain but not in others  

A translation invariant measure satisfies:
\begin{equation}
    M_\xi(f(\xi) + a) = M_\xi(f(\xi)) + a
    \label{eq:basic_trans_invar}
\end{equation}
if $a \in \mathbb{R}$ and $f(\xi) \in \mathscr{F}$. In a DO context, this property ensures that offsetting the cost function will not change the shape of the optimal cost trajectory. 

The positive homogeneity property is given by:
\begin{equation}
    M_\xi(\tau f(\xi)) = \tau M_\xi(f(\xi))
    \label{sec:basic_pos_homo}
\end{equation}
if $\tau > 0$ and $f(\xi) \in \mathscr{F}$. In the context of DO, a positive homogeneous  measure provides the property that uniformly scaling the cost by $\tau$ will not affect the shape of the optimal cost trajectory. 

The analysis of the stochastic risk measures featured in Section \ref{sec:measures} is well-established in the SO literature and Table \ref{tab:random_properties} provides a summary of these \citep{artzner1999coherent, ruszczynski2006optimization}.

\begin{table}[!htb]
    \centering
    \begin{tabular}{|c | c  c  c  c|}
        \hline
        $M_\xi$ & \eqref{eq:basic_convexity} & \eqref{eq:basic_monotonicity} & \eqref{eq:basic_trans_invar} & \eqref{sec:basic_pos_homo} \\ \hline
        $\mathbb{E}_\xi$ & Yes & Yes & Yes & Yes \\
        $\mathbb{E}$-$\mathbb{V}_\xi$ & Yes & No & Yes & No \\
        $Q_\xi$ & No & Yes & Yes & Yes \\
        CVaR$_\xi$ & Yes & Yes & Yes & Yes \\
        $\tilde{D}_\xi$ & Yes & Yes & Yes & Yes \\ \hline
    \end{tabular}
    \caption{A summary of the properties satisfied by certain measures $M_\xi$. Note that $\tilde{D}_\xi$ only satisfies \eqref{sec:basic_pos_homo} if $g(\cdot)$ is positive homogeneous.}
    \label{tab:random_properties}
\end{table}

Since the time-valued measures are special cases of the SO counterparts, they inherit the properties of Table \ref{tab:random_properties}. This illustrates how one can transfer rich theory from SO (with respect to these measure operators) to a DO context.  The time-valued expectation measure is a coherent risk measure; this might explain why this has been the classical measure used in DO. It is particularly important to observe that convexity ensures that the use of this measure yields a convex objective if the cost function is convex. Interestingly, the convexity of the objective has key implications for establishing stability conditions (e.g., closed-loop stability of MPC) \citep{rawlings2017model}. From this, we observe that other non-convex measures such as $\mathbb{E}$-$\mathbb{V}_t$ and $Q_t$ may not yield stability. On the other hand, CVaR$_t$ and $\tilde{D}_t$ are convex and thus might inherit stability properties (this is an interesting topic of future work). 

\subsection{Modeling in \texttt{InfiniteOpt.jl}} \label{sec:infiniteopt}
The unifying abstraction for infinite-optimization (which facilitated the connection between Problems \eqref{eq:gen_form} and \eqref{eq:2stage}) is implemented in a \texttt{Julia} package called \texttt{InfiniteOpt.jl} \citep{pulsipher2021unifying}. This enables us to intuitively model continuous DO problems (e.g., Problem \eqref{eq:gen_form}) following a simple symbolic syntax (see Code Snippet \ref{code:infiniteopt}). It is measure-centric and readily enables to quickly implement new candidate measure operators (such as those proposed in Section \ref{sec:measures}). Moreover, the unifying abstraction behind \texttt{InfiniteOpt.jl} facilities the incorporation of random constructs and/or PDE constraints. By default, these models are solved via direct transcription, but other solution methodologies can be implemented.

 
\section{Case Study} \label{sec:cases}

We compare the time-valued measures proposed in Section \ref{sec:measures} in the context of optimal control. We adapt the pandemic control problem that seeks to choose an isolation policy to control the spread of a contagion that minimally impacts the economic impact (imposed by mandated isolation). We model the spread of the disease via the SEIR model which defines the populations of individuals susceptible to infection $y_{s}: \mathcal{D}_{t} \rightarrow[0,1]$, exposed individuals that are not yet infectious $y_{e}: \mathcal{D}_{t} \rightarrow[0,1]$, infectious individuals $y_{i}: \mathcal{D}_{t} \rightarrow[0,1]$, and recovered individuals $y_{r}: \mathcal{D}_{t} \rightarrow[0,1]$ (considered immune to future infection). Moreover, these satisfy $y_{s}(t)+y_{e}(t)+y_{i}(t)+y_{r}(t)=1$. Thus, our state variables are comprised of $y_{s}(t)$, $y_{e}(t)$, $y_{i}(t)$, and $y_{r}(t)$. Moreover, we exhibit control by imposing an isolation policy $y_u(t) \in [0,\overline{y_u}] \subseteq [0, 1]$ that entails the separation of susceptible and exposed individuals ($y_u(t) = 0$ denotes no separation and $y_u(t) = 1$ denotes complete separation). The formulation seeks to minimize the isolation policy function $y_u(t)$ while enforcing that the amount of infectious individuals $y_i(t)$ remains below $\overline{y_i} \in (0, 1]$:
\begin{equation}
    \begin{aligned}
    && \min &&& M_t y_u(t) \\
    && \text{s.t.} &&& \frac{d y_{s}(t)}{d t}=\left(y_{u}(t)-1\right) \beta y_{s}(t) y_{i}(t), \ t \in \mathcal{D}_{t} \\
    &&&&& \frac{d y_{e}(t)}{d t}=\left(1-y_{u}(t)\right) \beta y_{s}(t) y_{i}(t)-\xi y_{e}(t), \ t \in \mathcal{D}_{t} \\
    &&&&& \frac{d y_{i}(t)}{d t}=\xi y_{e}(t)-\gamma y_{i}(t), \ t \in \mathcal{D}_{t}  \\
    &&&&& \frac{d y_{r}(t)}{d t}=\gamma y_{i}(t), \ t \in \mathcal{D}_{t} \\
    &&&&& y_{s}(0)=s_{0}, y_{e}(0)=e_{0}, y_{i}(0)=i_{0}, y_{r}(0)=r_{0}  \\
    &&&&&  y_{i}(t) \leq \overline{y_i}, \ t \in \mathcal{D}_{t} \\
    &&&&& y_{u}(t) \in\left[0, \overline{y_{u}}\right], \ t \in \mathcal{D}_{t}
    \end{aligned}
    \label{eq:pandemic_form}
\end{equation}
where $s_{0}, e_{0}, i_{0}, r_{0} \in [0, 1]$ are initial conditions and $\beta, \gamma, \xi \in \mathbb{R}$ are the rates of infection, recovery, and incubation, respectively, which are specific to the disease in question. Here we specify initial conditions at $s_{0} = .9999$, $e_{0}= 10^{-5}$, and $i_{0}=r_{0}= 0$. The disease parameters are taken to be $\beta = 0.727$, $\gamma=0.303$, and $\xi = 0.3$. We choose limits $\overline{y_u} = 0.8$ and $\overline{y_i} = 0.02$. Finally, we set $\mathcal{D}_t = [0, 200]$. We model Formulation \eqref{eq:pandemic_form} in \texttt{InfiniteOpt.jl} and use backward finite-difference to evaluate the derivatives using 101 discretization points. Code Snippet \ref{code:infiniteopt} highlights the compact syntax required to model Formulation \eqref{eq:pandemic_form} in \texttt{InfiniteOpt.jl} under these conditions. For the objective $M_t y_u(t)$, we consider the following time-valued measures: $\int_{t \in \mathcal{D}_t} y_u(t) dt$, $\mathbb{E}_t[y_u(t)]$, $\mathbb{E}$-$\mathbb{V}_t(y_u(t); \lambda)$, and CVaR$_t(y_u(t); \alpha)$. We also investigate the implications of using the uniform time-valued pdf in \eqref{eq:uniform_pdf} against the exponential pdf of \eqref{eq:exp_pdf}.

\begin{figure*}[!htp]
    \centering
    \begin{tabular}{cc}
      \includegraphics[width=.45\textwidth]{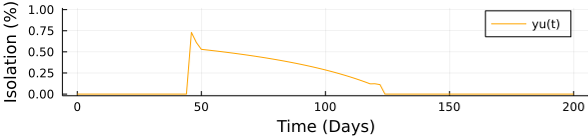}&
      \includegraphics[width=.45\textwidth]{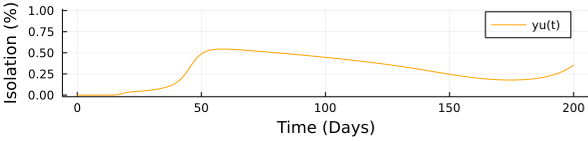}\\
      \includegraphics[width=.45\textwidth]{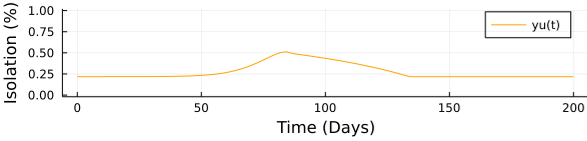}&
      \includegraphics[width=.45\textwidth]{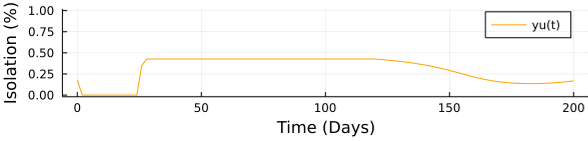}\\
    \end{tabular}
    \vspace{-0.1in}
    \caption{The optimal policy trajectories $y_u(t)$. Top Left: $\mathbb{E}_t$ with uniform pdf. Top Right: $\mathbb{E}_t$ with exponential pdf ($\gamma = 1$). Bottom Left: $\mathbb{E}$-$\mathbb{V}_t$ with $\lambda = 8$ and uniform pdf. Bottom Right: CVaR$_t$ with $\alpha = 0.9$ and uniform pdf.}
    \label{fig:results}
\end{figure*}

\begin{minipage}[!htb]{\linewidth}
\begin{scriptsize}
\lstset{language=Julia,breaklines = true}
\begin{lstlisting}[label = {code:infiniteopt},caption = Formulation \eqref{eq:pandemic_form} implemented in \texttt{InfiniteOpt.jl}.,captionpos=b]
using InfiniteOpt, Ipopt

# Set the parameters
γ, β, ξ = 0.303, 0.727, 0.3
s0, e0, i0, r0 = 1 - 1e-5, 1e-5, 0, 0

# Define the model
m = InfiniteModel(Ipopt.Optimizer)

# Define the time parameter
@infinite_parameter(m, t ∈ [0, 200], num_supports = 101)

# Add the variables
@variable(m, ys, Infinite(t))
@variable(m, ye, Infinite(t))
@variable(m, yi ≤ 0.02, Infinite(t))
@variable(m, yr, Infinite(t))
@variable(m, 0 ≤ yu ≤ 0.8, Infinite(t))

# Set the time expectation objective
@objective(m, Min, ��(yu, t))

# Define the SEIR equations
@constraint(m, ∂(ys, t) == -(1 - yu) * β * ys * yi)
@constraint(m, ∂(ye, t) == (1 - yu) * β * ys * yi - ξ * ye)
@constraint(m, ∂(yi, t) == ξ * ye - γ * yi)
@constraint(m, ∂(yr, t) == γ * yi)
@constraint(m, ys(0) == ys0)
@constraint(m, ye(0) == ye0)
@constraint(m, yi(0) == yi0
@constraint(m, yr(0) == yr0)

# Solve the model and retrieve results
optimize!(m)
u_opt = value(yu)
ts = value(t)
\end{lstlisting}
\end{scriptsize}
\end{minipage}

Figure \ref{fig:results} summarizes the results; the expectation $\mathbb{E}_t[y_u(t)]$ with the uniform pdf shapes $y_u(t)$ identically to the classical integral measure. However, we are able to place increased emphasis on the early time regime when we use the exponential pdf defined in \eqref{eq:exp_pdf} in combination with $\mathbb{E}_t[y_u(t)]$. In comparison to the other expectation measure, the exponentially weighted counterpart exhibits a policy trajectory that is significantly reduced at early times while later times lead to increased isolation requirements. This highlights how the choice of pdf $p_t(t)$ enhances the flexibility of our proposed measures in accordance with the requirements of the problem. The optimal policy we obtain with the mean-variance measure $\mathbb{E}$-$\mathbb{V}_t(y_u(t); 8)$ demonstrates how placing increased priority on minimizing the variance of the cost function induces the trajectory to be increasingly smoothed (i.e., cost fluctuations are damped). This comes at the trade-off (controlled via specification of $\lambda$) of increasing the mean isolation policy, but helps to derive a more consistent policy. For this application, a smoother policy would likely be preferred since rapid policy changes can be highly disruptive and can lead to public dissatisfaction. Finally, in contrast to the mean-variance (which equally penalizes positive and negative cost deviations from the mean),  the CVaR measure only penalizes the high cost deviations that surpass the threshold determined by the $\alpha$-quantile. In Figure \ref{fig:results} we see that CVaR$_t(y_u(t; 0.9)$ flattens the peak isolation policy values observed with the standard integral/expectation policy. This hedging against high costs also induces a more substantial response at later times which results in a larger cumulative cost. Thus, we observe a trade-off (controlled via $\alpha$) between penalizing cost peaks and minimizing the total cumulative cost. 

\section{Conclusions} \label{sec:conclusion}
In this work, we have shown that risk measures used in SO can be transferred (along with their mathematical properties) to DO. This enables a new class of DO formulations that shape time trajectories in interesting and useful ways. The transfer of insights across the SO and DO disciplines is facilitated by a unifying infinite-dimensional abstraction. In future work, it will be interesting to investigate the analogy between SO and DO further in transferring more amenable measure operators and establishing their properties in a DO context (e.g., dominance and stability). Moreover, the establishment of time-valued pdfs provides a foundation from which the utility of transferring distributionally robust measure functions to DO can be investigated; which would potentially allow us to consider multiple weighting functions over the time horizon.

\bibliography{references}

\begin{thebibliography}{12}
\providecommand{\natexlab}[1]{#1}
\providecommand{\url}[1]{\texttt{#1}}
\providecommand{\urlprefix}{URL }
\expandafter\ifx\csname urlstyle\endcsname\relax
  \providecommand{\doi}[1]{doi:\discretionary{}{}{}#1}\else
  \providecommand{\doi}{doi:\discretionary{}{}{}\begingroup
  \urlstyle{rm}\Url}\fi

\bibitem[{Artzner et~al.(1999)Artzner, Delbaen, Eber, and
  Heath}]{artzner1999coherent}
Artzner, P., Delbaen, F., Eber, J.M., and Heath, D. (1999).
\newblock Coherent measures of risk.
\newblock \emph{Mathematical finance}, 9(3), 203--228.

\bibitem[{Fulga(2016)}]{fulga2016portfolio}
Fulga, C. (2016).
\newblock Portfolio optimization with disutility-based risk measure.
\newblock \emph{European Journal of Operational Research}, 251(2), 541--553.

\bibitem[{Krokhmal et~al.(2013)Krokhmal, Zabarankin, and
  Uryasev}]{krokhmal2013modeling}
Krokhmal, P., Zabarankin, M., and Uryasev, S. (2013).
\newblock Modeling and optimization of risk.
\newblock \emph{Handbook of the fundamentals of financial decision making: Part
  II}, 555--600.

\bibitem[{Leland(1999)}]{leland1999beyond}
Leland, H.E. (1999).
\newblock Beyond mean--variance: Performance measurement in a nonsymmetrical
  world (corrected).
\newblock \emph{Financial analysts journal}, 55(1), 27--36.

\bibitem[{Pulsipher et~al.(In Press 2021)Pulsipher, Zhang, Hongisto, and
  Zavala}]{pulsipher2021unifying}
Pulsipher, J.L., Zhang, W., Hongisto, T.J., and Zavala, V.M. (In Press 2021).
\newblock A unifying modeling abstraction for infinite-dimensional
  optimization.
\newblock \emph{Computers \& Chemical Engineering}.

\bibitem[{Rawlings et~al.(2017)Rawlings, Mayne, and Diehl}]{rawlings2017model}
Rawlings, J.B., Mayne, D.Q., and Diehl, M. (2017).
\newblock \emph{Model predictive control: theory, computation, and design},
  volume~2.
\newblock Nob Hill Publishing Madison, WI.

\bibitem[{Renteria et~al.(2018)Renteria, Cao, Dowling, and
  Zavala}]{renteria2018optimal}
Renteria, J.A., Cao, Y., Dowling, A.W., and Zavala, V.M. (2018).
\newblock Optimal pid controller tuning using stochastic programming
  techniques.
\newblock \emph{AIChE Journal}, 64(8), 2997--3010.

\bibitem[{Risbeck and Rawlings(2019)}]{risbeck2019economic}
Risbeck, M.J. and Rawlings, J.B. (2019).
\newblock Economic model predictive control for time-varying cost and peak
  demand charge optimization.
\newblock \emph{IEEE Transactions on Automatic Control}, 65(7), 2957--2968.

\bibitem[{Rockafellar et~al.(2000)Rockafellar, Uryasev
  et~al.}]{rockafellar2000optimization}
Rockafellar, R.T., Uryasev, S., et~al. (2000).
\newblock Optimization of conditional value-at-risk.
\newblock \emph{Journal of risk}, 2, 21--42.

\bibitem[{Ruszczy{\'n}ski and Shapiro(2006)}]{ruszczynski2006optimization}
Ruszczy{\'n}ski, A. and Shapiro, A. (2006).
\newblock Optimization of risk measures.
\newblock In \emph{Probabilistic and randomized methods for design under
  uncertainty}, 119--157. Springer.

\bibitem[{Sarykalin et~al.(2008)Sarykalin, Serraino, and
  Uryasev}]{sarykalin2008value}
Sarykalin, S., Serraino, G., and Uryasev, S. (2008).
\newblock Value-at-risk vs. conditional value-at-risk in risk management and
  optimization.
\newblock In \emph{State-of-the-art decision-making tools in the
  information-intensive age}, 270--294. Informs.

\bibitem[{Shin and Zavala(2020)}]{shin2020diffusing}
Shin, S. and Zavala, V.M. (2020).
\newblock Diffusing-horizon model predictive control.
\newblock \emph{arXiv preprint arXiv:2002.08556}.

\end{thebibliography}
\end{document}